\documentclass[11pt,a4paper]{article}

\usepackage{epsf,epsfig,amsfonts,amsgen,amsmath,amstext,amsbsy,amsopn,amsthm,color
}
\usepackage{pgf,tikz}
\usepackage{enumitem}
\usepackage{url}
\usepackage{mathtools}
\usepackage{authblk}
\usepackage{amssymb}
\usepackage{dsfont}
\usepackage[backref=page]{hyperref}

\usepackage{fullpage}

\newtheorem{definition}{Definition} [section]
\newtheorem{theorem}[definition]{Theorem}

\newtheorem{corollary}[definition]{Corollary}

\newtheorem{claim}[definition]{Claim}

\newcommand{\graph}{H}

\begin{document}
\title{\bf\Large Sabotaging Mantel's Theorem}
\date{\today}
\author[1]{Natalie Behague\thanks{Research supported by the European Research Council (ERC) under the European Union Horizon 2020 research and innovation programme (grant agreement No. 947978).}\thanks{Email: \texttt{natalie.behague@warwick.ac.uk}}}
\author[2]{Debsoumya Chakraborti\protect\footnotemark[1]\thanks{Email: \texttt{cdebsoumya@gmail.com}}
}
\author[3]{Xizhi Liu\thanks{Research supported by ERC Advanced Grant 101020255 and the Excellent Young Talents Program (Overseas) of the National Natural Science Foundation of China. Email: \texttt{liuxizhi@ustc.edu.cn}}}
\affil[1]{Mathematics Institute,
            University of Warwick, 
            Coventry, UK}
\affil[2]{Department of Mathematics, University College London, London, UK}
\affil[3]{School of Mathematical Sciences, 
            University of Science and Technology of China, 
            Hefei, China}
\maketitle
\begin{abstract}
One of the earliest results in extremal graph theory, Mantel's theorem, states that the maximum number of edges in a triangle-free graph $G$ on $n$ vertices is $\lfloor n^2/4 \rfloor$. 
We investigate how this extremal bound is affected when $G$ is additionally required to contain a prescribed graph $\graph$ as a subgraph. 
We establish general upper and lower bounds for this problem, which are tight in the exponent for random triangle-free graphs and graphs generated by the triangle-free process, when the size of $\graph$ lies within certain ranges. 
\end{abstract}
\section{Introduction}\label{SEC:Introduction}

Mantel's theorem~\cite{Mantel07}, one of the earliest results in extremal graph theory, states that the number of edges in every triangle-free graph on $n$ vertices is at most $\lfloor n^2/4 \rfloor$, with equality attained  uniquely by the balanced complete bipartite graph $K_{\lfloor n/2\rfloor,\lceil n/2\rceil}$. 

We consider the following variation of this classical result: 
Given a vertex set $V$ of size $n$ and a prescribed triangle-free graph $\graph$ on $V$, what is the maximum number of edges in another triangle-free graph $G$ on $V$ that contains $\graph$ as a subgraph?
In other words, we are interested in the function 
\begin{align*}
    \mathrm{ex}_{\graph}(n,K_{3})
    \coloneqq \max\left\{e(G) \colon \text{$G \subseteq \binom{V}{2}$ is $K_{3}$-free and $\graph \subseteq G$}\right\}.
\end{align*}
Here, we identify a graph $G$ with its edge set and use $e(G)$ to denote the number of edges in it. 
It is clear that if $\graph$ is a subgraph of $K_{\lfloor n/2\rfloor,\lceil n/2\rceil}$, then $\mathrm{ex}_{\graph}(n,K_{3})$ equals $\mathrm{ex}(n,K_{3})$, and is strictly smaller otherwise (due to the fact that $K_{\lfloor n/2\rfloor,\lceil n/2\rceil}$ is uniquely extremal). Informally speaking, we will focus in this paper on the case where $H$ is `large'.

One of the simplest examples is perhaps when $\graph = S_{m}$, the star with $m$ edges.
Note that for $m > \lceil n/2 \rceil$, the graph $S_{m}$ is not a subgraph of $K_{\lfloor n/2\rfloor,\lceil n/2\rceil}$, so it follows from the discussion above that $\mathrm{ex}_{S_{m}}(n,K_{3}) < \mathrm{ex}(n,K_{3})$. 
The exact value of $\mathrm{ex}_{S_{m}}(n,K_{3})$ follows from a theorem\footnote{In fact, they showed that the same conclusion holds for large $n$ when $K_3$ is replaced by any cycle of odd length.} of Balister--Bollob{\'a}s--Riordan--Schelp in~\cite{BBRS03}, which shows that
\begin{align*}
    \mathrm{ex}_{S_{m}}(n,K_{3}) = m(n-m)
    \quad\text{for every $m > \lceil n/2 \rceil$.}
\end{align*}
The value of $\mathrm{ex}_{\graph}(n,K_{3})$ for a general triangle-free graph $\graph$ does not appear to have been systematically studied in the literature before.  
In this short note, we aim to establish general lower and upper bounds for $\mathrm{ex}_{\graph}(n,K_{3})$.
Before stating the main result, we introduce several related parameters. 

Given a graph $\graph$, let $d(\graph)$ denote its \textbf{average degree} and let $\Delta(\graph)$ denote its \textbf{maximum degree}. Let $\mathrm{N}(S_2, \graph)$ denote the number of (unlabelled) copies of $S_2$ in $\graph$, noting that
\begin{align}\label{ineq:N(S_2) vs Delta}
    e(\graph) + \mathrm{N}(S_2, \graph)
     = \frac{1}{2} \sum_{v\in V(\graph)} d_{\graph}(v) +  \sum_{v\in V(\graph)} \binom{d_{\graph}(v)}{2} 
    & = \frac{1}{2} \sum_{v\in V(\graph)} d_{\graph}(v)^2 \nonumber \\
    &\le \frac{1}{2} \sum_{v\in V(\graph)} d_{\graph}(v) \Delta(\graph)
    = e(\graph)  \Delta(\graph).
\end{align}
For a graph $\graph$ satisfying $e(\graph) + \mathrm{N}(S_2, \graph) =  \frac{1}{2} \sum_{v\in V(\graph)} d_{\graph}(v)^2  < \lfloor n^2/4 \rfloor$, we define the following function: 
\begin{align*}
    \gamma(\graph) 
    \coloneqq \frac{n(n-2)}{\lfloor n^2/4 \rfloor - e(\graph) - \mathrm{N}(S_2, \graph)} = \frac{n(n-2)}{\lfloor n^2/4 \rfloor - \frac{1}{2} \sum_{v\in V(\graph)} d_{\graph}(v)^2}. 
\end{align*}
The function $\psi$ may seem mysterious to the reader, and its exact statement arises simply as an artifact of our proof.

Recall that a vertex set $I \subseteq V(\graph)$ is an \textbf{independent set} of $\graph$ if it does not contain any edge of $\graph$. 
The \textbf{independence number} $\alpha(\graph)$ of $\graph$ is the size of the largest independent set in $\graph$. 
Let $\psi \colon [0,\infty) \to (0,1]$ be the decreasing continuous function defined by 
\begin{align*}
    \psi(d)
    \coloneqq 
    \begin{cases}
        \frac{ d \ln d - d + 1  }{(d-1)^2}, & \quad\text{if}\quad d \not\in \{0,1\} \\
        1,  & \quad\text{if}\quad d = 0, \\
        \frac{1}{2}, & \quad\text{if}\quad d = 1.
    \end{cases}
\end{align*}
The function $\psi$ comes from the following classical result of Shearer, which gives a lower bound on the independence number of triangle-free graphs.
\begin{theorem}[\cite{She83}]\label{THM:Shearer83}
    Suppose that $G$ is a triangle-free graph on $n$ vertices with average degree $d$. 
    Then $\alpha(G) \ge n \cdot \psi(d)$. 
\end{theorem}

Our main result is the following theorem. We use the notation $[n] \coloneq \{1,2,\ldots,n\}$. 
\begin{theorem}\label{THM:ex-B-general-bounds}
    Let $\graph$ be a $K_{3}$-free graph on vertex set $[n]$. 
    Then the following statements hold. 
    \begin{enumerate}[label=(\roman*)]
        \item\label{THM:ex-B-general-bounds-a} We have $\mathrm{ex}_{\graph}(n, K_{3}) \le {n \alpha(\graph)}/{2}$. 
        \item\label{THM:ex-B-general-bounds-b} Suppose that $ \frac{1}{2} \sum_{v\in V(\graph)} d_{\graph}(v)^2  < \lfloor n^2/4 \rfloor$. Then 
        \begin{align*}
            \mathrm{ex}_{\graph}(n, K_{3})
            \ge \left(\lfloor n^2/4 \rfloor -  \frac{1}{2} \sum_{v\in V(\graph)} d_{\graph}(v)^2 \right) \cdot \psi\big( \gamma(\graph) d(\graph)\big).
        \end{align*}
    \end{enumerate}
\end{theorem}

For positive real numbers $\beta_1$ and $\beta_2$, we say a graph $\graph$ on $[n]$ is \textbf{$(\beta_1, \beta_2)$-constrained} (or simply \textbf{constrained} when the constants $\beta_1$ and $\beta_2$ are not essential) if 
\begin{align*}
    \alpha(\graph) \le  \frac{\beta_1 n \ln d(\graph)}{d(\graph)} 
    \quad\text{and}\quad 
    e(\graph) \Delta(\graph) \le \left(\frac{1}{4} - \beta_2 \right) n^2.
\end{align*}
Note that for the first condition in the definition to hold, we implicitly require that $d(\graph) > 1$.

The following corollary is an immediate consequence of Theorem~\ref{THM:ex-B-general-bounds} and inequality~\eqref{ineq:N(S_2) vs Delta}. 

\begin{corollary}\label{CORO:ramsey-graph}
    Let $\beta_1 > 0$ and $\beta_2 > 0$ be two constants. 
    Suppose that $\graph$ is a $(\beta_1, \beta_2)$-constrained triangle-free graph on $[n]$ with average degree $d$. 
    Then 
    \begin{align*}
        \mathrm{ex}_{\graph}(n, K_{3}) 
        = \Theta_{n}\left(\frac{n^2 \ln  d}{d}\right).
    \end{align*}
\end{corollary}
A proof of Theorem~\ref{THM:ex-B-general-bounds} as well as the calculations giving Corollary~\ref{CORO:ramsey-graph} are provided in full in Section 2.

\bigskip

The definition of constrained graphs is motivated by properties of the Erd\H{o}s--R\'{e}nyi random graph $G(n,d/n)$, where each pair of vertices forms an edge independently with probability $d/n$. We consider this model for $d = d(n)$ be a fixed function of $n$ and we say  that an event occurs with high probability (\textit{w.h.p.}~for short) if it occurs with probability tending to $1$ as $n$ tends to infinity.
It is a well-known result in random graph theory (see e.g.~\cite{Fri90} and~{\cite[Theorem~3.4]{FK16book}}) that \textit{w.h.p.}, 
\begin{align*}
    \alpha(G(n,d/n)) 
    \le \frac{(2+o(1)) n \ln d}{d}
    \quad\text{and}\quad
    \Delta(G(n,d/n))
    \le (1+o(1)) \max\left\{d,~\frac{\ln n}{\ln\ln n}\right\}.
\end{align*}
Thus $G(n,d/n)$ is $(2+o(1), \varepsilon)$-constrained \textit{w.h.p.}~whenever $d \le \sqrt{(1/2 - \varepsilon) n}$.
Unfortunately, $G(n,d/n)$ is not triangle-free \textit{w.h.p.}~when $d \gg 1$ (see e.g.~{\cite[Theorem~1.12]{FK16book}}).
Below, we examine three random models that are triangle-free and share certain properties with the Erd\H{o}s--R\'{e}nyi random graph.

Let $\mathcal{T}(n,d)$ denote the family of all triangle-free graphs with vertex set $[n]$ and average degree $d$. 
Let $T(n,d)$ be a graph chosen uniformly at random from $\mathcal{T}(n,d)$. 
It was shown in~{\cite[Lemma~3]{OPT01}} that there exists $c > 0$ such that for $d \in [4, cn^{1/4}\sqrt{\ln  n} ]$, \textit{w.h.p.}, 
\begin{align}\label{equ:Tnd-independence-number}
    \alpha(T(n,d)) \le \frac{4 n \ln  d}{d}. 
\end{align}
By combining the Chernoff inequality (see e.g.~{\cite[Theorem~27.6]{FK16book}}) with~{\cite[Lemma~4]{OPT01}}, it is straightforward to verify that there exist constants $c', C> 0$ such that for every $d \in [4,c' n^{1/2}]$, \textit{w.h.p.},
\begin{align}\label{equ:Tnd-max-degree}
    \Delta(T(n,d))
    \le C d^3. 
\end{align}
Thus we obtain the following corollary for the case $\graph = T(n,d)$. 

\begin{corollary}\label{CORO:random-K3-free-graph}
    There exists a constant $c_{\ref{CORO:random-K3-free-graph}} > 0$ such that for every $d \in [4, c_{\ref{CORO:random-K3-free-graph}} n^{1/4}]$, \textit{w.h.p.}
    \begin{align*}
        \mathrm{ex}_{T(n,d)}(n, K_{3}) 
        = \Theta_{n}\left(\frac{n^2 \ln  d}{d}\right).
    \end{align*}
\end{corollary}

Another important model of random triangle-free graphs, which plays a significant role in constructions related to Ramsey theory for triangles (see, e.g.,~\cite{Boh09,FGGM20,BK21}), is the triangle-free process.  This process begins with the empty graph $G(0)$ on $n$ vertices. 
For each $i \ge 1$, the graph $G(i)$ is obtained by adding an edge $e_i$ to $G(i-1)$, where $e_i$ is chosen uniformly at random from the set of non-edges of $G(i-1)$ that do not form a triangle in $G(i)$.

It was shown in~\cite{Boh09,FGGM20,BK21} that the triangle-free process terminates \textit{w.h.p.}~after $\Theta(n^{3/2} \sqrt{\ln n})$ steps. 
Heuristically, the graph $G(i)$ closely resembles the Erd\H{o}s--R\'{e}nyi random graph $G(n,p)$ with $p = i/\binom{n}{2}$, except for the triangle-freeness. 
In particular, the maximum degree of $G(i)$ was shown in e.g.~\cite[Section~3.2]{BK21} to have \textit{w.h.p.}~the same order as those $G(n,p)$ with $p = i/\binom{n}{2}$, when $i$ is not too large. The same coupling argument can be repeated to show that this also holds for the independence number\footnote{Personal communication with Peter Keevash.}.

Specifically, there exists a constant $c> 0$ such that for $i \le c n^{3/2}$, \textit{w.h.p.}, and with $p = i/\binom{n}{2}$, 
\begin{align*}
    \Delta(G(i)) = O\left(\Delta(G(n,p))\right)
    \quad\text{and}\quad 
    \alpha(G(i)) = O\left(\alpha(G(n,p))\right).
\end{align*}
These properties, together with Corollary~\ref{CORO:ramsey-graph}, yield the following result for the case $\graph = G(i)$. 
\begin{corollary}\label{CORO:K3-free-process}
    There exists a constant $c_{\ref{CORO:K3-free-process}} > 0$ such that for every $i \in [n, c_{\ref{CORO:K3-free-process}} n^{3/2}]$ and with $d = 2i/n$,  \textit{w.h.p.}
    \begin{align*}
        \mathrm{ex}_{G(i)}(n, K_{3}) 
        = \Theta_{n}\left(\frac{n^2 \ln d}{d}\right).
    \end{align*}
\end{corollary}

\section{Proof of Theorem~\ref{THM:ex-B-general-bounds}}\label{SEC:Proof}
We prove Theorem~\ref{THM:ex-B-general-bounds} and Corollary~\ref{CORO:ramsey-graph} in this section.

\begin{proof}[Proof of Theorem~\ref{THM:ex-B-general-bounds}]
    Fix a $K_{3}$-free graph $\graph$ on $[n]$.     
    We say a graph $G$ on $[n]$ is $\graph$-admissible if it is $K_{3}$-free and contains $\graph$ as a subgraph. 
    
    First, we prove Theorem~\ref{THM:ex-B-general-bounds}~\ref{THM:ex-B-general-bounds-a}. 
    Suppose that $G$ is a $\graph$-admissible graph. 
    Let $v \in [n]$ be a vertex of maximum degree in $G$. 
    Since $G$ is $K_{3}$-free, the neighborhood $N_{G}(v)$ of $v$ in $G$ must be an independent set in $G$, and in particular, is independent in $\graph$.
    It follows that $\Delta(G) \le \alpha(\graph)$, and hence, 
    \begin{align*}
        e(G)
        \le \frac{n \cdot \Delta(G)}{2}  
        \le \frac{n \cdot \alpha(\graph)}{2},
    \end{align*}
    which proves Theorem~\ref{THM:ex-B-general-bounds}~\ref{THM:ex-B-general-bounds-a}. 

    Next, we prove Theorem~\ref{THM:ex-B-general-bounds}~\ref{THM:ex-B-general-bounds-b}. 
    Consider the auxiliary $3$-graph $\mathcal{H}$ whose vertex set is $\binom{[n]}{2}$ (i.e. the edge set of the complete graph $K_n$), where $\{e_1, e_2, e_3\} \subseteq \binom{[n]}{2}$ forms an edge in $\mathcal{H}$ iff $\{e_1, e_2, e_3\}$ spans a copy of $K_3$ in $K_n$.
    Define 
    \begin{align*}
        \mathcal{B}_{1}
        & \coloneqq \left\{e_1 \in \binom{[n]}{2} \colon \text{there exist $e_2, e_3 \in \graph$ such that $\{e_1, e_2, e_3\} \in \mathcal{H}$ }\right\} \quad\text{and}\quad \\[0.5em]
        \mathcal{B}_{2}
        & \coloneqq \left\{ \{e_1, e_2\} \subseteq \binom{[n]}{2} \colon \text{there exists $e_3 \in \graph$ such that $\{e_1, e_2, e_3\} \in \mathcal{H}$ }\right\},
    \end{align*}
    where we think of $\mathcal{B}_2$ as a graph on vertex set $\binom{[n]}{2}$.

    Let $\mathcal{B}_{3}$ denote the induced subgraph of $\mathcal{H}$ on $\binom{[n]}{2} \setminus \graph$.

    We say that a set $I \subseteq V(\mathcal{H})$ is independent in $\mathcal{H}$ if no edge of $\mathcal{H}$ is entirely contained in $I$. Note that every $K_3$-free graph on $[n]$ corresponds to an independent set in $\mathcal{H}$, and vice versa.
    In particular, $\graph$ is an independent set in $\mathcal{H}$, and a $\graph$-admissible graph is simply an independent set in $\mathcal{H}$ that contains $\graph$. 

    The following claim follows directly from the definitions. 

    \begin{claim}\label{CLAIM:G-B-independent}
        A graph $G$ is $\graph$-admissible iff $G = \graph \cup I$ for some $I \subseteq \binom{[n]}{2} \setminus \graph$ satisfying all of the following conditions: 
        \begin{enumerate}[label=(\roman*)]
            \item\label{CLAIM:G-B-independent-1} $I \cap \mathcal{B}_{1} = \emptyset$, 
            \item\label{CLAIM:G-B-independent-2} $I$ is an independent set in $\mathcal{B}_{2}$, and 
            \item\label{CLAIM:G-B-independent-3} $I$ is an independent set in $\mathcal{B}_{3}$. 
        \end{enumerate}
    \end{claim}

    Therefore, to prove the lower bound for $\mathrm{ex}_{\graph}(n,K_{3})$, it suffices to provide a lower bound for $I \subseteq \binom{[n]}{2} \setminus \graph$ satisfying Claim~\ref{CLAIM:G-B-independent}~\ref{CLAIM:G-B-independent-1},~\ref{CLAIM:G-B-independent-2}, and~\ref{CLAIM:G-B-independent-3}.

    \begin{claim}\label{CLAIM:B1-size-upper-bound}
        We have $e(\mathcal{B}_{1}) \le \mathrm{N}(S_{2}, \graph)$.
    \end{claim}
    \begin{proof}[Proof of Claim~\ref{CLAIM:B1-size-upper-bound}]
        It follows from the definition of $\mathcal{B}_{1}$ that a pair $\{u,v\} \in \binom{[n]}{2}$ is contained in $\mathcal{B}_{1}$ iff there exists a vertex $w \in [n] \setminus \{u,v\}$ such that $\{w,u\} \in \graph$ and $\{w,v\} \in \graph$. In other words, $\{u,v\}$ is contained in the neighborhood $N_{\graph}(w)$ of some vertex $w$. 
        Therefore, 
        \begin{align*}
            e(\mathcal{B}_{1})
            \le \sum_{w\in [n]} \binom{d_{\graph}(w)}{2}
            = \mathrm{N}(S_2, \graph),
        \end{align*}
        which proves Claim~\ref{CLAIM:B1-size-upper-bound}. 
    \end{proof}

    \begin{claim}\label{CLAIM:B2-size-upper-bound}
        We have $e(\mathcal{B}_{2}) = e(\graph)(n-2)$.
    \end{claim}
    \begin{proof}[Proof of Claim~\ref{CLAIM:B2-size-upper-bound}]
        It follows from the definition that $\{e_1, e_2\} \subseteq \binom{[n]}{2}$ is an edge in $\mathcal{B}_{2}$ iff there exists three vertices $u, v, w \in [n]$ such that $e_1 = \{u,v\}$, $e_2 = \{u,w\}$, and $\{v,w\} \in \graph$. 
        Therefore, each edge in $\graph$ contributes exactly $n-2$ edges to $\mathcal{B}_{2}$. 
        It follows that $e(\mathcal{B}_{2}) = e(\graph)(n-2)$. 
    \end{proof}

    \begin{claim}\label{CLAIM:B2-K3-free}
        The graph $\mathcal{B}_{2}$ is $K_{3}$-free. 
    \end{claim}
    \begin{proof}[Proof of Claim~\ref{CLAIM:B2-K3-free}]
        Suppose to the contrary that there exist three pairs $e_1, e_2, e_3 \in \binom{[n]}{2}$ that span a copy of $K_{3}$ in $\mathcal{B}_{2}$, i.e. $\big\{\{e_1, e_2\}, \{e_1, e_3\}, \{e_2, e_3\} \big\} \subseteq \mathcal{B}_{2}$. 
        By the definition of $\mathcal{B}_{2}$, for each pair of indices $\{i,j\} \in \binom{[3]}{2}$, there exists $e_{ij} \in \graph$ such that $\{ e_i, e_j, e_{ij} \}$ forms a triangle in $K_{n}$, i.e. $|e_i \cap e_j| = 1$ and $e_{ij} = e_i \triangle e_j$, the symmetric difference of $e_i$ and $e_j$. 

        Since $|e_1 \cap e_2| = |e_1 \cap e_3| = |e_2 \cap e_3| = 1$, the set $\{e_1, e_2, e_3\}$ must span either a $3$-edge star or a triangle in $K_{n}$. 
        In both cases, the corresponding set $\{e_{12}, e_{13}, e_{23}\} \subseteq \graph$ forms a triangle, contradicting the $K_{3}$-freeness of $\graph$. 
        Therefore, $\mathcal{B}_{2}$ is $K_{3}$-free.
    \end{proof}

    We are now ready to construct the set $I \subseteq \binom{[n]}{2} \setminus \graph$ that satisfies Claim~\ref{CLAIM:G-B-independent}~\ref{CLAIM:G-B-independent-1},~\ref{CLAIM:G-B-independent-2}, and~\ref{CLAIM:G-B-independent-3}.  
    We begin with an independent set $S \subseteq \binom{[n]}{2}$ in $\mathcal{H}$ of size $\lfloor n^2/4 \rfloor$. Note that as this corresponds to a $K_{3}$-free graph on $[n]$ of size $\lfloor n^2/4 \rfloor$, it is uniquely realized by choosing a balanced complete bipartite graph on $[n]$.
    
    Define $S' \coloneqq S \setminus (\graph \cup \mathcal{B}_{1})$, and let $I$ be a maximum independent set of the induced subgraph $\mathcal{B}_{2}[S']$. Note that $I$ satisfies Claim~\ref{CLAIM:G-B-independent}~\ref{CLAIM:G-B-independent-1},~\ref{CLAIM:G-B-independent-2}, and~\ref{CLAIM:G-B-independent-3}.  

    \begin{claim}\label{CLAIM:I2-lower-bound}
        We have $|I| \ge \left(\lfloor n^2/4 \rfloor - e(\graph) - \mathrm{N}(S_2, \graph)\right) \cdot \psi\big( \gamma(\graph) d(\graph)\big)$.
    \end{claim}
    \begin{proof}[Proof of Claim~\ref{CLAIM:I2-lower-bound}]
        Let $d$ denote the average degree of the induced subgraph $\mathcal{B}_{2}[S']$. 
        It follows from Claims~\ref{CLAIM:B2-size-upper-bound} and~\ref{CLAIM:B1-size-upper-bound} that 
        \begin{align*}
            d 
            = \frac{2 e(\mathcal{B}_{2}[S'])}{|S'|}
            \le \frac{2 e(\mathcal{B}_{2})}{|S'|}
            & \le \frac{2 e(\graph) (n-2)}{\lfloor n^2/4 \rfloor - e(\graph) - e(\mathcal{B}_{1})} \\
            & \le \frac{2 e(\graph) (n-2)}{\lfloor n^2/4 \rfloor - e(\graph) - \mathrm{N}(S_2, \graph)}
            = \gamma(\graph)  d(\graph).
        \end{align*}
        Since, by Claim~\ref{CLAIM:B2-K3-free}, $\mathcal{B}_{2}[S']$ is $K_{3}$-free, it follows from Theorem~\ref{THM:Shearer83} that 
        \begin{align*}
            |I| 
            = \alpha(\mathcal{B}_{2}[S'])
            \ge |S'| \cdot \psi(d)
            \ge \left(\left\lfloor \frac{n^2}{4} \right\rfloor - e(\graph) - \mathrm{N}(S_2, \graph)\right) \cdot \psi\big( \gamma(\graph)  d(\graph)\big),
        \end{align*}
        which proves Claim~\ref{CLAIM:I2-lower-bound}. 
    \end{proof}

    Claim~\ref{CLAIM:I2-lower-bound} completes the proof of Theorem~\ref{THM:ex-B-general-bounds}~\ref{THM:ex-B-general-bounds-b}.
\end{proof}

Now we will provide a full proof of the Corollary.
\begin{proof}[Proof of Corollary~\ref{CORO:ramsey-graph}]
    Let  $\graph$ be a $(\beta_1, \beta_2)$-constrained triangle-free graph on $[n]$ with average degree $d$. By Theorem~\ref{THM:ex-B-general-bounds}~\eqref{THM:ex-B-general-bounds-a}, we have \[\mathrm{ex}_{\graph}(n, K_{3}) \le {n \alpha(\graph)}/{2} \le \frac{\beta_1 }{2}\cdot \frac{n^2 \ln d}{d}.\] Applying inequality~\eqref{ineq:N(S_2) vs Delta}, we see that $\frac{1}{2} \sum_{v\in V(\graph)} d_{\graph}(v)^2 \le e(\graph) \Delta(\graph) \le \left(\frac{1}{4} - \beta_2 \right) n^2.$ 
    As a result, clearly $\beta_2 < 1/4$ and \[\gamma(H) = \frac{n(n-2)}{\lfloor n^2/4\rfloor - \frac{1}{2} \sum_{v\in V(\graph)} d_{\graph}(v)^2} \le \frac{n^2}{\lfloor n^2/4\rfloor - \left(\frac{1}{4} - \beta_2 \right) n^2} \le \frac{2}{\beta_2}.\]
   As $\psi$ is decreasing, this implies that 
    \[\psi\big( \gamma(\graph) d\big) 
    \ge   \frac{\frac{2}{\beta_2} d\ln\left({\frac{2}{\beta_2}d}\right) - \frac{2}{\beta_2}d + 1}{\left(\frac{2}{\beta_2}d - 1\right)^2} 
    \ge   \frac{\frac{2}{\beta_2} d\ln({d}) + \frac{2}{\beta_2}d\left(\ln(\frac{2}{\beta_2}) - 1\right)}{\left(\frac{2}{\beta_2}d\right)^2} 
    \ge \frac{\beta_2}{2} \frac{\ln{d}}{d}.\]
    Therefore, by Theorem~\ref{THM:ex-B-general-bounds}~\eqref{THM:ex-B-general-bounds-b}, we have            
    \begin{align*} \mathrm{ex}_{\graph}(n, K_{3})
            &\ge \left(\lfloor n^2/4 \rfloor -  \frac{1}{2} \sum_{v\in V(\graph)} d_{\graph}(v)^2 \right) \cdot \psi\big( \gamma(\graph) d\big) \\
            &\ge  \left(\lfloor n^2/4 \rfloor -  \left(\frac{1}{4} - \beta_2 \right) n^2 \right) \cdot \frac{\beta_2}{2} \frac{\ln{d}}{d} \\
            &\ge \frac{\beta_2^2}{4} \cdot \frac{n^2\ln{d}}{d}.
            \end{align*}
    This completes the proof.
\end{proof}
\section{Concluding remarks}\label{SEC:Remark}
\begin{enumerate}
    \item The study of the function $\mathrm{ex}_{\graph}(n,K_{3})$ was partially motivated by results in~\cite[Section~5]{CLSW25}, where the authors needed to upper bound the number of edges in graphs satisfying certain properties and containing a specific subgraph. Both of these problems fall under the scope of the following more general meta-question:
What happens to an extremal problem when we require the presence of a prescribed subgraph/subset?
We hope that our results could inspire further research in this direction. 
\item Determining the asymptotic behavior of $\mathrm{ex}_{T(n,d)}(n, K_{3})$ and $\mathrm{ex}_{G(i)}(n, K_{3})$ for all feasible values of $d$ and $i$ seems to be an interesting open problem. 
It is well known in random graph theory that when $d = o(1)$, the components of the random graph $T(n,d)$ (which is essentially equivalent to the  Erd\H{o}s--R\'{e}nyi random graph $G(n,d/n)$ in this regime) are, \textit{w.h.p.}, trees of size $o(\ln  n)$ (see~{\cite[Theorems~2.1~and~2.9]{FK16book}}).
As a result, $T(n,d)$ admits a bipartition of its vertex set in which the part sizes differ by $o(\ln  n)$. Therefore, \textit{w.h.p.}, 
\begin{align*}
    \mathrm{ex}_{T(n,d)}(n,K_{3}) = (1/4 - o(1))n^2
    \quad\text{if}\quad d = o(1).
\end{align*}
For $d \ge \left( \sqrt{3}/2 + o(1) \right) \sqrt{n \ln n}$, a classical theorem of Osthus--Pr{\"o}mel--Taraz~\cite{OPT03} shows that, \textit{w.h.p.}, $T(n,d)$ is bipartite. Moreover, their proof implies that the bipartition has nearly equal part sizes, that is, the part sizes differ by $o(n)$.
Therefore, \textit{w.h.p.}~
\begin{align*}
    \mathrm{ex}_{T(n,d)}(n,K_{3}) = (1/4 - o(1))n^2
    \quad\text{if}\quad d \ge \left( \sqrt{3}/2 + o(1) \right) \sqrt{n \ln n}.
\end{align*}
The cases not covered by Corollary~\ref{CORO:random-K3-free-graph} and the discussion above remain open for $\mathrm{ex}_{T(n,d)}(n, K_{3})$.  
It is worth noting that recent work by Jenssen--Perkins--Potukuchi~\cite{JPP23} appears to be helpful for the case when $d$ is near the threshold $\left( \sqrt{3}/2 + o(1) \right) \sqrt{n \ln n}$. 
Similarly, the range not addressed by Corollary~\ref{CORO:K3-free-process} remains open for $\mathrm{ex}_{G(i)}(n, K_{3})$.

\item  In~\cite{GL23}, Guo--Warnke introduced a refined alteration model in which one constructs an $F$-free graph by deleting all edges lying in copies of $F$ in $G(n,p)$, rather than removing just a single edge from each copy as in classical alteration arguments. They showed that, for appropriate values of $p$, this full-deletion approach does not substantially affect key parameters such as the independence number (see Theorem~1 and Lemma~6 therein), and often leads to cleaner analyses in applications such as online Ramsey games.
One can apply Theorem~\ref{THM:ex-B-general-bounds} to their model and obtain straightforwardly an analogous result to that of Corollary~\ref{CORO:K3-free-process}. 

Warnke also pointed out to us that the structural properties needed for Corollary~\ref{CORO:K3-free-process} follow from standard random-graph estimates and do not require more sophisticated analysis of the triangle-free process.

\item A natural extension of Theorem~\ref{THM:ex-B-general-bounds} is to replace $K_3$ with a general family of graphs $\mathcal{F}$. 
Given a family $\mathcal{F}$ of graphs, a graph $G$ is \textbf{$\mathcal{F}$-free} if it does not contain any member of $\mathcal{F}$ as a subgraph.  
Analogously, given an $\mathcal{F}$-free graph $\graph$ on $[n]$, one may ask for the value of the following function: 
\begin{align*}
    \mathrm{ex}_{\graph}(n,\mathcal{F})
    \coloneqq \left\{e(G) \colon \text{$G \subseteq \binom{[n]}{2}$ is $\mathcal{F}$-free and $\graph \subseteq G$}\right\}.
\end{align*}

\item Note that Corollaries~\ref{CORO:random-K3-free-graph} and~\ref{CORO:K3-free-process} can be viewed as results concerning $\mathrm{ex}_{\graph}(n,K_{3})$ in the average case. 
One can also consider the worst-case scenario. 
Let $m, n \ge 1$ be integers. Define 
\begin{align*}
    \mathrm{ex}_{m}(n,\mathcal{F}) 
    = \min\left\{ \mathrm{ex}_{\graph}(n,\mathcal{F}) \colon \text{$\graph$ is $\mathcal{F}$-free and $e(\graph) \le m$} \right\}.
\end{align*}
The function $\mathrm{ex}_{m}(n,\mathcal{F})$ is closely related to a classical problem introduced by Erd{\H o}s--Hajnal--Moon~\cite{EHM64}. 
Given a family $\mathcal{F}$ of graphs, a graph $G$ is \textbf{$\mathcal{F}$-saturated} if $G$ is $\mathcal{F}$-free but the addition of any new edge to $G$ creates a subgraph that belongs to $\mathcal{F}$. 
The \textbf{saturation number} $\mathrm{sat}(n,\mathcal{F})$ is defined as 
\begin{align*}
    \mathrm{sat}(n,\mathcal{F})
    \coloneqq \min\left\{e(G) \colon \text{$v(G) = n$ and $G$ is $\mathcal{F}$-saturated}\right\}.
\end{align*}
Note that for $m \ge \mathrm{sat}(n,\mathcal{F})$, the function $\mathrm{ex}_{m}(n,\mathcal{F})$ reduces to $\mathrm{sat}(n,\mathcal{F})$.  

For $m$ close to $\mathrm{sat}(n,K_3)$ (i.e. $n-1$), it appears that $S_{m}$ is the worst graph, that is, $\mathrm{ex}_{m}(n,K_3) = \mathrm{ex}_{S_m}(n,K_3) = m(n-m)$. However, for smaller values of $m$, the disjoint union of $5$-cycles is worse than $S_m$. 
We leave it as an open problem to determine the behavior of $\mathrm{ex}_{m}(n,\mathcal{F})$, in particular $\mathrm{ex}_{m}(n,K_3)$, as $m$ increases from $1$ to $\mathrm{sat}(n,\mathcal{F})$. 

\item The definition of the function $\mathrm{ex}_{\graph}(n,\mathcal{F})$ can be extended even to the case when $\graph$ is not necessarily $\mathcal{F}$-free: 
Given a graph $\graph$ on a vertex set $V$ of size $n$ and a family of graphs $\mathcal{F}$, what is the maximum number of edges in a graph $G$ on $V$ that contains $\graph$ as a subgraph and satisfies $\mathrm{N}(F,G) = \mathrm{N}(F,\graph)$ for every $F\in \mathcal{F}$? (For graphs $G_1$ and $G_2$, the notation $\mathrm{N}(G_1,G_2)$ denotes the number of copies of $G_1$ in $G_2$.) 
Here, we are interested in the function 
\begin{align*}
    \mathrm{ex}_{\graph}(n,\mathcal{F})
    \coloneqq \max\left\{e(G) \colon \text{$G \subseteq \binom{V}{2}$, $\graph \subseteq G$, and $\mathrm{N}(F,G) = \mathrm{N}(F,\graph)$ for every $F\in \mathcal{F}$}\right\}.
\end{align*}
For the general function $\mathrm{ex}_{\graph}(n,K_3)$ when $\graph$ is not necessarily triangle-free, Corollary~\ref{CORO:ramsey-graph} continues to hold, which then applies for the random graph models, where triangle-freeness is not required. A slight modification of our arguments can be used to prove Corollary~\ref{CORO:ramsey-graph} in this more general setting. We omit the details and leave these general problems to future investigation. 

\end{enumerate}

\section*{Acknowledgements}
We would like to thank Peter Keevash for clarifying certain properties of the triangle-free process. 
We would also like to thank Lutz Warnke for drawing our attention to~\cite{GL23} and pointing out its connection to our work.

\bibliographystyle{alpha}
\bibliography{sabotageturan}
\end{document}